%&amstex
\input amsppt.sty
\hsize=12.5cm
\vsize=540pt
\NoRunningHeads
\TagsOnRight
%-----------------------------------------
  \define\CC{\Bbb C}
  \define\di{\diamond}
\redefine\epsilon{\varepsilon}
  \define\intt{\operatorname{int}}
  \define\NN{\Bbb N}
  \define\one{\bold1}
  \define\RR{\Bbb R}
  \define\skipaline{\vskip12pt}
  \define\sprod#1#2{\langle#1,#2\rangle}
  \define\too{\longrightarrow}
  \define\two{\bold2}
  \define\ZZ{\Bbb Z}
%-----------------------------------------
\topmatter

\title
On $\bold N$-circled $\bold{H^\infty}$-domains of holomorphy
\endtitle

\author
Marek Jarnicki (Krak\'ow), Peter Pflug (Oldenburg)
\endauthor

\address
Uniwersytet Jagiello\'nski \newline
Instytut Matematyki \newline
30-059 Krak\'ow, Reymonta 4, Poland\newline
{\it E-mail: }{\rm jarnicki{\@}im.uj.edu.pl}
\endaddress

\address
Carl von Ossietzky Universit\"at Oldenburg
\newline Fachbereich Mathematik\newline
Postfach 2503\newline
D-26111 Oldenburg, Germany\newline
{\it E-mail: }{\rm pflugvec{\@}dosuni1.rz.uni-osnabrueck.de}\newline
{\it E-mail: }{\rm pflug{\@}hrz2.pcnet.uni-oldenburg.de}
\endaddress

\thanks
Research supported by KBN Grant No 2 PO3A 060 08
and by Volkswagen Stiftung Az\. I/71~062
\endthanks

\endsubjclass

\endkeywords

\abstract
We present various characterizations of $n$-circled domains of holomorphy
$G\subset\CC^n$ with respect to some subspaces of $\Cal H^\infty(G)$.
\endabstract

\endtopmatter

\document

%==========================================

\head Introduction \endhead
We say that a domain $G\subset\CC^n$ is {\it $n$-circled\/} if
$(e^{i\theta_1}z_1,\dots,e^{i\theta_n}z_n)\in G$ for arbitrary
$(z_1,\dots,z_n)\in G$ and $(\theta_1,\dots,\theta_n)\in\RR^n$.

Put $\log G:=\{(x_1,\dots,x_n)\in\RR^n\: (e^{x_1},\dots,e^{x_n})\in G\}$.

If $X\subset\RR^n$ is a convex domain, then $\Cal E(X)$ denotes
the largest vector subspace $F\subset\RR^n$ such that $X+F=X$.

A vector subspace $F\subset\RR^n$ is said to be of {\it rational type\/}
if $F$ spanned by $F\cap\ZZ^n$.

Let $L_h^2(G):=\Cal O(G)\cap L^2(G)$.
\skipaline
The following results are known (cf\. \cite{Jar-Pfl 1}).

\proclaim{Proposition 1} Let $G\subset\CC^n$ be an $n$-circled domain of
holomorphy. Then the following conditions are equivalent:

\noindent{\rm (i)} $G$ is fat (i.e\. $G=\intt\overline G$)
and the space $\Cal E(\log G)$ is of rational type;

\noindent{\rm (ii)} there exist $A\subset\ZZ^n$ and
$(c_\alpha)_{\alpha\in A}\subset(0,+\infty)$ such that
$$
G=\intt\bigcap_{\alpha\in A}\{|z^\alpha|<c_\alpha\};
$$

\noindent{\rm (iii)} $G$ is an $\Cal H^\infty(G)$-domain of holomorphy.
\endproclaim

\proclaim{Proposition 2} Let $G\varsubsetneq\CC^n$ be a fat $n$-circled
domain of holomorphy. Then the following conditions are equivalent:

\noindent{\rm (i)} $\Cal E(\log G)=\{0\}$;

\noindent{\rm (ii)} $L_h^2(G)\neq\{0\}$;

\noindent{\rm (iii)} $G$ is an $L_h^2(G)$-domain of holomorphy.
\endproclaim

Define
$$
L_h^{p,k}(G):=\{f\in\Cal O(G)\: \forall_{|\nu|\leq k}\:
\partial^\nu f\in L^p(G)\},\quad p\in[1,+\infty],\;k\in\ZZ_+,
$$
where $\partial^\nu:=\frac{\partial^{|\nu|}}{\partial z_1^{\nu_1}\dots
\partial z_n^{\nu_n}}$, $\nu=(\nu_1,\dots,\nu_n)\in(\ZZ_+)^n$.
Note that the space $L_h^{p,k}(G)$ endowed with the norm
$\|f\|_{L_h^{p,k}(G)}:=\sum_{|\nu|\leq k}\|\partial^\nu f\|_{L^p(G)}$
is a Banach space. Put
$$
L_h^p(G):=L_h^{p,0}(G)=\Cal O(G)\cap L^p(G),\quad
\Cal H^{\infty,k}(G):=L_h^{\infty,k}(G).
$$
Note that $L_h^\infty(G)=\Cal H^{\infty,0}(G)=\Cal H^\infty(G)$. Let
$$
\gather
L_h^{\di,k}(G):=\bigcap_{p\in[1,+\infty]}L_h^{p,k}(G),\\
\Cal A^k(G):=\{f\in\Cal O(G)\: \forall_{|\nu|\leq k}\;
\exists_{f_\nu\in\Cal C(\overline G)}\: f_\nu=\partial^\nu f \text{ in } G\},
\quad k\in\ZZ_+,\\
\Cal A^\infty(G):=\bigcap_{k\in\ZZ_+}\Cal A^k(G).
\endgather
$$

\proclaim{Remark}\;\rm (a) If the volume of $G$ is finite, then
$\Cal H^{\infty,k}(G)\subset L_h^{\di,k}(G)$.

\noindent (b) If $G$ is bounded, then
$\Cal A^k(G)\subset\Cal H^{\infty,k}(G)$.
\endproclaim

\skipaline
The aim of this paper is to generalize Propositions 1, 2. The starting point
of these investigations was our trial to understand the general situation
behind the last example in \S\;I.5 in \cite{Sib}. We will prove the
following results.

\proclaim{Proposition 3} Let $G\subset\CC^n$ be an $n$-circled domain
of holomorphy. Then the following conditions are equivalent:

\noindent{\rm (i)} $G$ is fat and $\Cal E(\log G)=\{0\}$;

\noindent{\rm (ii)} $G$ is fat and there exists
$p\in[1,+\infty)$ such that $L_h^p(G)\neq\{0\}$;

\noindent{\rm (iii)} $G\varsubsetneq\CC^n$ and for each $k\in\ZZ_+$
the domain $G$ is an $L_h^{\di,k}(G)\cap\Cal A^k(G)$-domain of holomorphy.
\endproclaim

Observe that for each $k\in\ZZ_+$ the space
$\Cal F_k(G):=L_h^{\di,k}(G)\cap\Cal A^k(G)$ endowed with the norms
$\|\;\|_{L_h^{p,k}(G)}$, $p\in[1,+\infty]$, is a Fr\'echet space.
Consequently, condition (iii) is equivalent to the following one.

\noindent(iii') {\sl For each $k\in\ZZ_+$ there exists a function
$f\in L_h^{\di,k}(G)\cap\Cal A^k(G)$ such that $G$ is the domain of
existence of $f$.}

Note that $\Cal F_k(G)$ is an algebra.

\proclaim{Corollary 4} Let $G\subset\CC^n$ be a bounded fat $n$-circled
domain of holomorphy. Then $G$ is an $\Cal A^k(G)$-domain
of holomorphy (for any $k\in\ZZ_+$).
\endproclaim

In the case where $G=\{|z_1|<|z_2|<1\}$ is the Hartogs triangle the above
result has been proved in \cite{Sib}.
\skipaline
For $\epsilon=(\epsilon_1,\dots,\epsilon_n)\in\{0,1\}^n$ put
$$
\align
V_\epsilon:=\{(z_1,\dots,z_n)\in\CC^n\:
&1^o\: z_j=0 \text{ for all $j$ such that }\epsilon_j=1,\\
&2^o\: z_j\neq0 \text{ for all $j$ such that }\epsilon_j=0\},
\endalign
$$

\vskip-12pt
$$
G_\epsilon:=\{(\lambda_1^{\epsilon_1}z_1,\dots,
\lambda_n^{\epsilon_n}z_n)\:
\lambda_1,\dots,\lambda_n\in\overline E,\; (z_1,\dots,z_n)\in G\},
$$
where $E$ denotes the unit disc.

\proclaim{Proposition 5} Let $G\subset\CC^n$ be an $n$-circled
$\Cal H^\infty(G)$-domain of holomorphy. Then the following conditions are
equivalent:

\noindent{\rm (i)} for any $\epsilon\in\{0,1\}^n$ if $V_\epsilon\cap\partial
G\neq\varnothing$, then $G_\epsilon\subset G$;

\noindent{\rm (ii)} $G$ is a $\Cal S(G)$-domain of holomorphy, where
$$
\multline
\Cal S(G):=\{f\in\Cal O(G)\: \sup_{K\cap G}|\partial^\nu f|<+\infty\\
\text{ for any compact $K\subset\overline G$ and for any $\nu\in(\ZZ_+)^n$}\};
\endmultline
$$

\noindent{\rm (iii)} $G$ is an $\Cal A^\infty(G)$-domain of holomorphy;

\noindent{\rm (iv)} $G$ is an
$\Cal H^\infty(G)\cap\Cal O(\overline G)$-domain of holomorphy.
\endproclaim

\proclaim{Corollary 6} Let $G\subset\CC^n$ be a bounded fat $n$-circled
domain of holomorphy. Then $G$ is an  $\Cal A^\infty(G)$-domain of holomorphy
iff $G$ satisfies condition (i) of Proposition 5.
\endproclaim

In the case where $G$ is the Hartogs triangle the above result has been
proved (by different methods) in \cite{Sib}.

\proclaim{Proposition 7} Let $G\subset\CC^n$ be an $n$-circled
domain of holomorphy. Then the following conditions are equivalent:

\noindent{\rm (i)} $G=D\times\CC^{n-m}$, where $D$ is a fat $m$-circled
domain of holomorphy with $\Cal E(\log D)=\{0\}$ and $0\leq m\leq n$;

\noindent{\rm (ii)} $G$ is an $\Cal H^{\infty,1}(G)$-domain of holomorphy;

\noindent{\rm (iii)} for any $k\in\ZZ_+$ the domain $G$ is an
$\Cal H^{\infty,k}(G)\cap\Cal A^k(G)$-domain of holomorphy.
\endproclaim

%===========================================

\head Proof of Proposition 3 \endhead

(iii) $\Longrightarrow$ (ii) follows from Proposition 1.
\skipaline
(ii) $\Longrightarrow$ (i). Let $f=\sum_{\nu\in\ZZ^n}a_\nu z^\nu
\in L_h^p(G)$, $f\not\equiv0$. Then
$$
\align
\int_G|a_\nu z^\nu|^pd\Lambda_{2n}(z)
&=(2\pi)^n\int_{|G|}\bigg|\frac1{(2\pi i)^n}\int_{\Sb|\zeta_j|=r_j\\
j=1,\dots,n\endSb}\frac{f(\zeta)}{\zeta^{\nu+\one}}
d\zeta\bigg|^pr^{p\nu+\one}d\Lambda_n(r)\\
&\leq(2\pi)^{n(1-p)}\int_{|G|}\bigg(\int_{[0,2\pi]^n}
|f(re^{i\theta)}|d\Lambda_n(\theta)\bigg)^pr^{\one}d\Lambda_n(r)\\
&\leq\int_{|G|}\int_{[0,2\pi]^n}|f(re^{i\theta})|^pd\Lambda_n(\theta)\;
r^{\one}d\Lambda_n(r)\\
&=\int_G|f|^pd\Lambda_{2n},
\endalign
$$
where $|G|:=\{(|z_1|,\dots,|z_n|)\: (z_1,\dots,z_n)\in G\}$ and
$\Lambda_n$ denotes Lebesgue measure in $\RR^n$, $\one:=(1,\dots,1)$.
Consequently, there exists $\nu_0\in\ZZ^n$ such that $z^{\nu_0}\in L_h^p(G)$.
\skipaline
Suppose that $F:=\Cal E(\log G)\neq\{0\}$. Let $k:=\dim F$ and let
$Y\subset F^\perp$ be a convex domain such that $\log G=Y+F$. We have:
$$
\align
\int_G|z^{\nu_0}|^pd\Lambda_{2n}(z)&
=(2\pi)^n\int_{\log G}e^{\sprod{x}{p\nu_0+\two}}d\Lambda_n(x)\\
&=\int_Ye^{\sprod{x'}{p\nu_0 +\two}}d\Lambda_{n-k}(x')
\int_Fe^{\sprod{x''}{p\nu_0+\two}}d\Lambda_k(x'')\\
&=+\infty,
\endalign
$$
where $\sprod{\;}{\;}$ is the Euclidean scalar product in $\RR^n$,
$\two:=2\cdot\one$. We have got a contradiction.
\skipaline
(i) $\Longrightarrow$ (iii). Fix $k\in\ZZ_+$ and put
$\Cal F_k(G):=L_h^{\di,k}(G)\cap\Cal A^k(G)$.

Suppose that there exist domains
$G_0$, $\widetilde G\subset\CC^n$ such that $\varnothing\neq G_0\subset G
\cap\widetilde G$, $\widetilde G\not\subset G$, and for each $f\in\Cal F_k(G)$
there exists $\widetilde f\in\Cal O(\widetilde G)$ with $\widetilde f=f$
on $G_0$.

Since $G$ is fat, we may assume that $\widetilde G\not\subset
\overline G$ and that $\widetilde G\cap W_0=\varnothing$,
where $W_0:=\{(z_1,\dots,z_n)\in\CC^n\: z_1\cdot\dots\cdot z_n=0\}$.
We may also assume that $\one\in\widetilde G\setminus\overline G$.
For, if $(a_1,\dots,a_n)\in\widetilde G\setminus\overline G$, then we
replace $G$ by $F(G)$, where $F(z_1,\dots,z_n):=(z_1/a_1,\dots,z_n/a_n)$.

Since $\Cal E(\log G)=\{0\}$ and $G$ is fat,
there exist $\RR$-linearly independent vectors
$\alpha_j=(\alpha_{j,1},\dots,\alpha_{j,n})\in\ZZ^n$, $j=1,\dots,n$, such that
$$
G\subset D:=\{z\in U\:|z^{\alpha_j}|<1\:j=1,\dots,n\},
$$
where 
$$
U:=\{(z_1,\dots,z_n)\in\CC^n\: \text{ if for some $j,\ell$ we have
$\alpha_{j,\ell}<0$, then } z_\ell\neq0\}
$$
(cf\. \cite{Jar-Pfl 1}).
Define $\Phi\:U\too\CC^n$, $\Phi(z):=(z^{\alpha_1},\dots,z^{\alpha_n})$.
Note that $\Phi$ is holomorphic, $D=\Phi^{-1}(E^n)$, and that
$\Phi$ is biholomorphic in a neighbourhood of $\one$. In particular, there
exists a point $b\in\widetilde G\setminus\overline G$ such that
$r_j:=|b^{\alpha_j}|>1$, $j=1,\dots,n$.

Let $\gamma\:[0,1]\too\widetilde G$
be an arc with $\gamma(0)\in G_0$, $\gamma(1)=b$. Put
$$
D_0:=\{z\in U\: |z^{\alpha_j}|<r_j,\; j=1,\dots,n\}
$$
and let $t_0:=\sup\{t>0\: \gamma([0,t))\subset D_0\}$.
Put $c:=\gamma(t_0)$. Then there exists $j_0\in\{1,\dots,n\}$ such that
$|c^{\alpha_{j_0}}|=r_{j_0}$. We may assume that $j_0=n$.

Put $\alpha:=\alpha_1+\dots+\alpha_n$. For $N\in\NN$ define
$$
f_N(z):=\frac{z^{N\alpha}}{z^{\alpha_n}-c^{\alpha_n}},
\quad z\in D_0.
$$
Obviously $f_N\in\Cal O(D_0)$. Note that $\gamma([0,t_0))\subset
\widetilde G\cap D_0$ and $\lim_{t\to t_0-}f_N(\gamma(t))=\infty$. In
particular, $f_N$ cannot be holomorphically continued to $\widetilde G$.

Consequently, to get a contradiction it suffices to show that
$f_N\in\Cal F_k(D)$ provided $N\gg1$.

For any $\sigma\in(\ZZ_+)^n$ with $|\sigma|\leq k$ we get
$$
\partial^\sigma f_N(z)=
-\sum_{\mu=0}^\infty\sigma!\binom{N\alpha+\mu\alpha_n}{\sigma}
\frac{z^{N\alpha+\mu\alpha_n-\sigma}}{d^{\mu+1}},\quad z\in D_0,
$$
where $d:=c^{\alpha_n}$ (recall that $|d|=r_n>1$). One can easily show that
$$
\bigg|\sigma!\binom{N\alpha+\mu\alpha_n}{\sigma}\bigg|\leq (P+Q\mu)^R,
\quad\mu\in\ZZ_+,
$$
where $P, Q, R\in\NN$ depend only on $k, \alpha_1,\dots,\alpha_n$, and $N$.
Since the series
$$
\sum_{\mu=0}^\infty\frac{(P+Q\mu)^R}{r_n^{\mu+1}}
$$
is convergent, we only need to find $N\in\NN$ such that for arbitrary
$|\sigma|\leq k$ and $p\in[1,+\infty]$ we get
$$
z^{N\alpha+\mu\alpha_n-\sigma}\in L_h^p(D)\cap\Cal C(\overline D),
\quad\|z^{N\alpha+\mu\alpha_n-\sigma}\|_{L^p(D)}\leq 1.
$$
Observe that $|z^{N\alpha+\mu\alpha_n-\sigma}|\leq|z^{N\alpha-\sigma}|$ on
$D$. Hence, it is enough to prove that there exists $N\in\NN$ such that
for arbitrary $|\sigma|\leq k$ and $p\in[1,+\infty]$ we have
$$
\|z^{N\alpha-\sigma}\|_{L^p(D)}\leq 1,\quad
\lim_{D\ni z\to z_0}z^{N\alpha-\sigma}=0,\quad z_0\in W_0\cap\partial D.
$$
\skipaline
Let $A:=[\alpha_{j,\ell}]_{j,\ell=1,\dots,n}$, $B:=A^{-1}$,
$T_j(x):=(xB)_j=\sum_{\ell=1}^nB_{\ell,j}x_\ell$, $j=1,\dots,n$,
$x=(x_1,\dots,x_n)\in\RR^n$.
\skipaline
For $p\in[1,+\infty)$ and $\nu\in\ZZ^n$ we have
$$
\align
\int_D|z^\nu|^pd\Lambda_{2n}(z)&=(2\pi)^n\int_{\log D}e^{\sprod{x}{p\nu+
\bold2}}d\Lambda_n(x)\\
&=(2\pi)^n\int_{\{\xi<0\}}e^{\sprod{B\xi}
{p\nu+\bold2}}|\det B|d\Lambda_n(\xi)\\
&=\frac{(2\pi)^n}{|\det A|T_1(p\nu+\two)\cdot\dots\cdot T_n(p\nu+\two)}
\endalign
$$
provided that $T_j(p\nu+\two)>0$, $j=1\dots,n$. In particular, if
$$
T_j(\nu)\geq\frac1p\bigg(\frac{2\pi}{|\det A|^{1/n}}-T_j(\two)\bigg),
\quad j=1,\dots,n,
$$
then $\|z^\nu\|_{L^p(D)}\leq 1$. Hence, if $\nu=N\alpha-\sigma$ and if
$$
\multline
N\geq N_0:=\max\{
T_j(\sigma)+\frac1p\bigg(\frac{2\pi}{|\det A|^{1/n}}-T_j(\two)\bigg)\:\\
 j=1,\dots,n,\;\sigma\in(\ZZ_+)^n,\;|\sigma|\leq k+1,\;p\in[1,+\infty)\},
\endmultline
$$
then $\|z^{N\alpha-\sigma}\|_{L^p(D)}\leq1$ for arbitrary $p\in[1,+\infty)$
and $|\sigma|\leq k$.
\skipaline
Moreover, $N_0\geq T_j(\sigma)$, $j=1,\dots,n$, and therefore
$N\alpha-\sigma\in\RR_+\alpha_1+\dots+\RR_+\alpha_n$, which shows that
$\|z^{N\alpha-\sigma}\|_{\Cal H^\infty(D)}\leq1$ for arbitrary $|\sigma|
\leq k$.
\skipaline
Finally, let $z_0=(z_{0,1},\dots,z_{0,n})\in W_0\cap\partial D$. Assume
that $z_{0,\ell_0}=0$. Let $(e_1,\dots,e_n)$ denotes the canonical basis of
$\RR^n$. Let $|\sigma|\leq k$.
Since $N_0\geq T_j(\sigma+e_{\ell_0})$, $j=1,\dots,n$, we get
$N\alpha-\sigma-e_{\ell_0}\in\RR_+\alpha_1+\dots+\RR_+\alpha_n$. This shows
that $|z^{N\alpha-\sigma}|\leq|z_{\ell_0}|$ in $D$ and therefore
$\lim_{D\ni z\to z_0}z^{N\alpha-\sigma}=0$.
\qed
%==========================================

\head Proof of Proposition 5 \endhead

The implications (iv) $\Longrightarrow$ (iii) $\Longrightarrow$ (ii) are
trivial.
\skipaline
(ii) $\Longrightarrow$ (i). Suppose that $V_\epsilon\cap\partial G\neq
\varnothing$ for an $\epsilon\in\{0,1\}^n$, $\epsilon\neq(0,\dots,0)$.
We may assume that $\epsilon_1=\dots\epsilon_s=1$,
$\epsilon_{s+1}=\dots=\epsilon_n=0$ for some $1\leq s\leq n$.
We will show that each function
$$
f=\sum_{\nu\in\ZZ^n}a_\nu(f) z^\nu\in \Cal S(G)
$$
extends holomorphically to $G_\epsilon$. Consequently, since
$G$ is an $\Cal S(G)$-domain of holomorphy, we will get
$G_\epsilon\subset G$.

Fix an $f\in\Cal S(G)$. It suffices to show that
$a_\mu(f)=0$ for all $\mu=(\mu_1,\dots,\mu_n)\in\ZZ^n$ such that
$\mu_j<0$ for some $j\in\{1,\dots,s\}$.
We may assume that $\mu_1\geq0, \dots \mu_t\geq0$, $\mu_{t+1}<0,\dots,
\mu_s<0$, where $0\leq t\leq s-1$. Define $\sigma:=(\mu_1,\dots,\mu_t,0,
\dots,0)\in(\ZZ_+)^n$. Let $r_0=(r_{0,1},\dots,r_{0,n})\in
V_\epsilon\cap\partial G\cap(\RR_+)^n$ and let
$$
\gather
U:=\{(r_1,\dots,r_n)\in G\: |r_j-r_{0,j}|<1,\; r_j>0,\; j=1,\dots,n\},\\
K:=\{(e^{i\theta_1}r_1,\dots,e^{i\theta_n}r_n)\: (r_1,\dots,r_n)\in
\overline U,\; (\theta_1,\dots,\theta_n)\in\RR^n\},\\
M:=\sup_{K\cap G}|\partial^\sigma f|.
\endgather
$$
Since
$$
\partial^\sigma f(z)=\sum_{\nu\in\ZZ^n}\sigma!\binom\nu\sigma a_\nu(f)
z^{\nu-\sigma},
$$
the Cauchy inequalities give:
$$
\sigma!|a_\mu(f)|\leq\frac{M}{r^{\mu-\sigma}}=
\frac{M}{r_{t+1}^{\mu_{t+1}}\cdot\dots\cdot r_n^{\mu_n}},
\quad r=(r_1,\dots,r_n)\in U.
$$
Now, letting $U\ni r\too r_0$, we conclude that $a_\mu(f)=0$.
\skipaline
(i) $\Longrightarrow$ (iii). We know that $G$ is of the form
$$
G=\intt\bigcap_{\alpha\in A}\{|z^\alpha|<c_\alpha\}\tag 1
$$
with $A\subset\ZZ^n$ and $z^\alpha\in\Cal H^\infty(G)$, $\alpha\in A$
(cf\. Proposition 1). We will prove that
$z^\alpha\in\Cal O(\overline G)$ for arbitrary $\alpha\in A$.

Fix an $\alpha=(\alpha_1,\dots,\alpha_n)\in A$.
We may assume that $\alpha_1<0,\dots,\alpha_s<0$, $\alpha_{s+1}\geq0,\dots
\alpha_n\geq0$ for some $1\leq s\leq n$. We have to prove that
$\partial G\cap\{z_1\cdot\dots\cdot z_s=0\}=\varnothing$. Suppose that
$\partial G\cap\{z_1\cdot\dots\cdot z_s=0\}\neq\varnothing$.
Then we may assume that $\partial G\cap V_\epsilon\neq\varnothing$ for
some $\epsilon\in\{0,1\}^n$ with $\epsilon_1=1$.
By (iii) $G_\epsilon\subset G$. In particular,
$$
|\lambda^{\alpha_1}||z^\alpha|<c_\alpha,
\quad \lambda\in\overline E,\; z\in G;
$$
contradiction.

Now, suppose that there exist domains $G_0$,
$\widetilde G\subset\CC^n$ such that $\varnothing\neq G_0\subset G\cap
\widetilde G$, $\widetilde G\not\subset G$, and for each $f\in\Cal H^\infty(G)
\cap\Cal O(\overline G)$ there exists $\widetilde f\in\Cal O(\widetilde G)$
with $\widetilde f=f$ on $G_0$. Since $G$ is fat, we may assume that
$\widetilde G\cap\{z_1\cdot\dots\cdot z_n=0\}=\varnothing$.

First observe that $\|\widetilde f\|_{\widetilde G}\leq\|f\|_G$. For, suppose
that $|\widetilde f(a)|>\|f\|_G$ for some $f\in\Cal H^\infty(G)
\cap\Cal O(\overline G)$ and $a\in\widetilde G$.
Then the function $g:=1/(f-\widetilde f(a))$ belongs to $\Cal H^\infty(G)
\cap\Cal O(\overline G)$ and cannot be holomorphically continued to
$\widetilde G$; contradiction.

Consequently, $|z^\alpha|<c_\alpha$ on $\widetilde G$ for any $\alpha\in A$.
Hence, by (1), $\widetilde G\subset G$; contradiction.
\qed

%==========================================

\head Proof of Proposition 7 \endhead

(i) $\Longrightarrow$ (iii) follows from Proposition 3.
(iii) $\Longrightarrow$ (ii) is trivial.
\skipaline
(ii) $\Longrightarrow$ (i). Let $F:=\Cal E(\log G)$, $m:=\dim F$. The cases
$m=0$  and $m=n$ are trivial. Assume $1\leq m\leq n-1$. One can easily
prove (cf\. \cite{Jar-Pfl 1}) that for any $f\in\Cal H^\infty(G)$ the Laurent
series of $f$ has the form
$$
\sum_{\nu\in F^\perp\cap\ZZ^n}a_\nu(f)z^\nu,\quad z\in G.\tag 2
$$
Now, suppose that $f\in\Cal H^{\infty,1}(G)$. By the above remark (applied
to $f$ and $\partial f/\partial z_j$ simultaneously) we see that
if $\nu=(\nu_1,\dots,\nu_n)\in
\ZZ^n$ is such that $\nu_j\neq 0$ and $a_\nu(f)\neq 0$, then $\nu, \nu-e_j
\in F^\perp$, and, consequently, $e_j\in F^\perp$
($(e_1,\dots,e_n)$ denotes the canonical basis in $\RR^n$).
Since $\dim F^\perp=n-m$, we may assume that $e_{s+1},\dots,e_n\notin
F^\perp$ for some $0\leq s\leq n-m$. Thus, the series (2) is independent
of the variables  $z_{s+1},\dots,z_n$ and, therefore, it is convergent in the
domain $D\times\CC^{n-s}$, where $D$ denotes the projection of $G$ on $\CC^s$.
Since $G$ is an $\Cal H^{\infty,1}(G)$-domain of holomorphy,
$G=D\times\CC^{n-s}$ and $D$ is an
$\Cal H^{\infty,1}(D)$-domain of holomorphy.
Moreover, $\Cal E(\log G)=\Cal E(\log D)\times\RR^{n-s}$. Hence $n-s=m$ and
$\Cal E(\log D)=\{0\}$.
\qed

In \cite{Sic 1,2} J\. Siciak characterized those balanced domains of
holomorphy $G\subset\CC^n$, which are $\Cal H^\infty(G)$-  (resp\.
$\Cal H^\infty(G)\cap\Cal A^\infty(G)$-) domains of holomorphy. Moreover,
it is known that any bounded balanced domain of holomorphy $G\subset\CC^n$ is
an $L_h^2(G)$-domain of holomorphy (cf\. \cite{Jar-Pfl 2}). A general
discussion for balanced domains of holomorphy (like the one for $n$-circled
domains) is still unknown.
%==========================================
\Refs
\widestnumber\key{Jar-Pfl 2}
{
%------------------------------------------
\ref
\key Jar-Pfl 1
\by M\. Jarnicki \& P\. Pflug
\paper Existence domains of holomorphic functions of restricted growth
\jour Trans\. Amer\. Math\. Soc\.
\yr 1987
\vol 304
\pages 385--404
\endref
%------------------------------------------
\ref
\key Jar-Pfl 2
\by M\. Jarnicki \& P\. Pflug
\paper On balanced $L^2$-domains of holomorphy
\jour Ann\. Pol\. Math\.
\yr 1995
\toappear
\endref
%-----------------------------
\ref
\key Sib
\by N\. Sibony
\paper Prolongement de fonctions holomorphes born\'ees et metrique
de Carath\'eo\-do\-ry
\jour Invent\. Math\.
\yr 1975
\vol 29
\pages 205--230
\endref
%-----------------------------
\ref
\key Sic 1
\by J\. Siciak
\paper Circled domains of holomorphy of type $\Cal H^\infty$
\jour Bull\. Soc\. Sci\. et de Lettres de \L \'od\'z
\yr 1984
\vol 34
\pages 1--20
\endref
%-----------------------------
\ref
\key Sic 2
\by J\. Siciak
\paper Balanced domains of holomorphy of type $H^\infty$
\jour Mat\. Vesnik
\yr 1985
\vol 37
\pages 134--144
\endref
%-----------------------------
}
\endRefs

\enddocument